\documentclass{article}
\usepackage{amsfonts}
\usepackage{amsmath}
\usepackage{graphicx}

\begin{document}

\title{Improving Accuracy of Goodness-of-fit Test.}
\author{Kris Duszak and Jan Vrbik \\
Brock University}
\maketitle

\begin{abstract}
It is well known that the approximate distribution of the usual test
statistic of a goodness-of-fit test is chi-square, with degrees of freedom
equal to the number of categories minus $1$ (assuming that no parameters are
to be estimated -- something we do throughout this article). Here we show how
to improve this approximation by including two correction terms, each of
them inversely proportional to the total number of observations.
\end{abstract}

\section{\textbf{Goodness-of-fit Test: A Brief Review}}

To test whether a random independent sample of size $n$ comes from a
specific distribution can be done by dividing all possible outcomes of the
corresponding random variable (say $U$) into $k$ distinct regions (called 
\textsc{categories}) so that these have similar probabilities of happening.
The sample of $n$ values of $U$ is then converted into the corresponding
observed frequencies, one for each category (we denote these $%
X_{1},X_{2},...X_{k}$), equivalent to sampling a multinomial distribution
with probabilities $p_{1},p_{2},...p_{k}$.(computed, for each category,
based on the original distribution). The new random variables $X_{i}$ have
expected values given by $n\cdot p_{i}$ (where $i$ goes from $1$ to $k$) and
variance-covariance matrix given by%
\begin{equation*}
n\cdot (\mathbb{P}-\mathbf{p\ p}^{T})
\end{equation*}%
where $\mathbf{p}$ is a column vector with $k$ elements (the individual $%
p_{i}$ probabilities), and $\mathbb{P}$ is similarly an $k$ $\times $ $k$ 
\emph{diagonal} matrix, with the same $p_{i}$ probabilities on its main
diagonal.

The usual test statistic is%
\begin{equation}
T=\sum_{i=1}^{k}\frac{(X_{i}-n\cdot p_{i})^{2}}{n\cdot p_{i}}\equiv
\sum_{i=1}^{k}Y_{i}^{2}  \label{Tstat}
\end{equation}%
where%
\begin{equation}
Y_{i}\equiv \frac{X_{i}-n\cdot p_{i}}{\sqrt{n\cdot p_{i}}}  \label{yi}
\end{equation}%
equivalent to (in its vector form)%
\begin{equation}
\mathbf{Y}=\frac{\mathbb{P}_{k}^{-1/2}(\mathbf{X}-n\cdot \mathbf{p})}{\sqrt{n%
}}  \label{Y}
\end{equation}%
where $\mathbf{X}$ is a column vector of the $X_{1},X_{2},...,X_{k}$
observations.

The $Y_{i}$'s have a mean of zero and their variance-covariance matrix is%
\begin{equation}
\mathbb{V}=\mathbb{P}^{-1/2}(\mathbb{P}-\mathbf{p\ p}^{T})\mathbb{P}^{-1/2}=%
\mathbb{I}-\mathbf{p}^{1/2}(\mathbf{p}^{1/2})^{T}  \label{one}
\end{equation}%
where $\mathbb{I}$ is the $k$ $\times $ $k$ unit matrix and $\mathbf{p}%
^{1/2} $ denotes a column vector with elements equal to $%
p_{1}^{1/2},p_{2}^{1/2},...p_{k}^{1/2}$. The matrix (\ref{one}) is \emph{%
idempotent}, since%
\begin{equation*}
\mathbf{p}^{1/2}(\mathbf{p}^{1/2})^{T}\mathbf{p}^{1/2}(\mathbf{p}^{1/2})^{T}=%
\mathbf{p}^{1/2}(\mathbf{p}^{1/2})^{T}
\end{equation*}%
and its \emph{trace} is $k-1$, since%
\begin{equation*}
\text{Tr}\left[ \mathbf{p}^{1/2}(\mathbf{p}^{1/2})^{T}\right] =\text{Tr}%
\left[ (\mathbf{p}^{1/2})^{T}\mathbf{p}^{1/2}\right] =\sum_{i=1}^{k}p_{i}=1.
\end{equation*}%
Because the $k$-dimensional distribution of (\ref{Y}) tends (as $%
n\rightarrow \infty $) to a Normal distribution with zero means and
variance-covariance matrix of (\ref{one}), (\ref{Tstat}) must similarly
converge to the $\chi _{k-1}^{2}$ distribution (assuming that $U$ does have
the hypothesized distribution). A substantial disagreement between the
observed frequencies $X_{i}$ and their expected values $n\cdot p_{i}$ will
be reflected by the test statistic $T$ exceeding the (right-hand-tail)
critical value of $\chi _{k-1}^{2},$ leading to a rejection of the null
hypothesis.

Since the sample size is always finite, the critical value (computed under
the assumption that $n\rightarrow \infty $) with have an error roughly
proportional to $\frac{1}{n}.$ To remove this error is an objective of this
article.

\section{$\frac{1}{n}$\textbf{\ proportional correction}}

A small modification of the results of \cite{vrbik} indicate that a
substantially better approximation (which removes the $\frac{1}{n}$
-proportional error) to the probability density function (PDF) of the
distribution of $T$ (under the null hypothesis) is%
\begin{eqnarray}
&&\chi _{k-1}^{2}(t)\cdot \left( 1+B\cdot (\frac{t^{2}}{(k-1)(k+1)}-\frac{2t%
}{k-1}+1)+\right.  \label{corr} \\
&&\left. C\cdot (\frac{t^{3}}{(k-1)(k+1)(k+3)}-\frac{3t^{2}}{(k-1)(k+1)}+%
\frac{3t}{k-1}-1)\right)  \notag
\end{eqnarray}%
where $\chi _{k-1}^{2}(t)$ is the PDF of the regular chi-square distribution
and%
\begin{equation}
B=\frac{1}{8}\sum_{i,j=1}^{k}\kappa _{i,i,j,j}  \label{B}
\end{equation}%
$\qquad \qquad \qquad \qquad \qquad \qquad \qquad \qquad \qquad \qquad
\qquad \qquad \qquad \qquad \qquad \qquad \qquad \ \ \ \qquad \qquad \qquad
\qquad \qquad \qquad \qquad \qquad \qquad \qquad \qquad \qquad $%
\begin{equation}
C=\frac{1}{8}\sum_{i,j,\ell =1}^{k}\kappa _{i,j,j}\kappa _{i,\ell ,\ell }+%
\frac{1}{12}\sum_{i,j,\ell =1}^{k}\kappa _{i,j,\ell }^{2}  \label{C}
\end{equation}%
where $\kappa _{i,j,\ell }$ and $\kappa _{i,j,\ell ,h,}$, are cumulants of
the (multivariate) $\mathbf{Y}$ distribution. They can be found easily,
based on the logarithm of the joint moment generating function of (\ref{yi}%
), namely%
\begin{equation*}
M=n\cdot \ln \left( \sum_{m=1}^{k}p_{m}\exp \left( \frac{t_{m}}{\sqrt{n\cdot
p_{m}}}\right) \right) -\sum_{m=1}^{k}t_{m}\sqrt{n\cdot p_{m}}
\end{equation*}%
by differentiating $M$ with respect to $t_{i},$ $t_{j}$ and $t_{\ell }$ to
get $\kappa _{i,j,\ell }$ (and the extra $t_{h}$ to get $\kappa _{i,j,\ell }$%
), followed by setting all $t_{m}=0$.

This yields%
\begin{eqnarray*}
\kappa _{i,i,i} &=&\frac{(1-p_{i})(1-2p_{i})}{\sqrt{n\cdot p_{i}}} \\
\kappa _{i,i,j} &=&-\frac{\sqrt{p_{j}}(1-2p_{i})}{\sqrt{n}} \\
\kappa _{i,j,\ell } &=&\frac{2\sqrt{p_{i}\cdot p_{j}\cdot p_{\ell }}}{\sqrt{n%
}}
\end{eqnarray*}
and%
\begin{eqnarray*}
\kappa _{i,i,i,i} &=&\frac{(1-p_{i})(1-6p_{i}+6p_{i}^{2})}{n\cdot p_{i}}=%
\frac{1}{n}\left( \frac{1}{p_{i}}-7+12p_{i}-6p_{i}^{2}\right) \\
\kappa _{i,i,j,j} &=&\frac{2p_{i}+2p_{j}-6p_{i}\cdot p_{j}-1}{n}.
\end{eqnarray*}

Using these formulas, we can proceed to compute%
\begin{eqnarray*}
&&B\overset{}{=}\frac{1}{8}\sum_{i=1}^{k}\kappa _{i,i,i,i}+\frac{1}{8}%
\sum_{i\neq j}^{k}\kappa _{i,i,j,j}\overset{}{=} \\
&&\frac{1}{8n}\left( Q-7k+12s_{1}-6(s_{1}^{2}-2s_{2})\overset{}{+}%
2(k-1)s_{1}+2(k-1)s_{1}-12s_{2}-k(k-1)\right)
\end{eqnarray*}%
where%
\begin{equation*}
Q\equiv \sum_{i=1}^{k}\frac{1}{p_{i}}
\end{equation*}%
and $s_{1}$ and $s_{2}$ are the first two elementary symmetric polynomials
in $p_{i},$ i.e.%
\begin{eqnarray*}
s_{1} &=&\sum_{i=1}^{k}p_{i} \\
s_{2} &=&\sum_{i<j}^{k}p_{i}\cdot p_{j}
\end{eqnarray*}%
(note that $\sum_{i=1}^{k}p_{i}^{2}=s_{1}^{2}-2s_{2}$). Realizing that $%
s_{1}=1,$ the expression for $B$ can be simplified to%
\begin{equation}
B=\frac{1}{8n}\left( Q-k^{2}-2k+2\right)  \label{BBres}.
\end{equation}%
When choosing the categories in a manner which makes all $p_{i}$ equal to $%
1/k$, the last expression reduces to%
\begin{equation*}
-\frac{k-1}{4n}
\end{equation*}

Similarly,%
\begin{eqnarray*}
&&C\overset{}{=}\frac{1}{8}\sum_{i=1}^{k}\kappa _{i,i,i}^{2}+\frac{1}{4}%
\sum_{i\neq j}^{k}\kappa _{i,i,i}\kappa _{i,j,j}+\frac{1}{8}\sum_{i\neq
j}^{k}\kappa _{i,j,j}^{2}+\frac{1}{8}\sum_{i\neq j\neq \ell }^{k}\kappa
_{i,j,j}\kappa _{i,\ell ,\ell } \\
&&+\frac{1}{12}\sum_{i=1}^{k}\kappa _{i,i,i}^{2}+\frac{1}{4}\sum_{i\neq
j}^{k}\kappa _{i,i,j}^{2}+\frac{1}{12}\sum_{i\neq j\neq \ell }^{k}\kappa
_{i,j,\ell }^{2} \\
&&\overset{}{=}\frac{5}{24n}\sum_{i=1}^{k}\frac{(1-p_{i})^{2}(1-2p_{i})^{2}}{%
p_{i}}-\frac{1}{4n}\sum_{i\neq j}^{k}(1-p_{i})(1-2p_{i})(1-2p_{j}) \\
&&+\frac{3}{8n}\sum_{i\neq j}^{k}p_{j}(1-2p_{i})^{2}+\frac{1}{8n}\sum_{i\neq
j\neq \ell }^{k}p_{i}(1-2p_{j})(1-2p_{\ell })+\frac{1}{3n}\sum_{i\neq j\neq
\ell }^{k}p_{i}p_{j}p_{\ell } \\
&&\overset{}{=}\frac{5}{24n}\sum_{i=1}^{k}\left( \frac{1}{p_{i}}%
-6+13p_{i}-12p_{i}^{2}+4p_{i}^{3}\right) \\
&&-\frac{1}{4n}\sum_{i=1}^{k}\left( k(1-3p_{i}+2p_{i}^{2})\overset{}{-}%
3+11p_{i}-12p_{i}^{2}+4p_{i}^{3}\right) \\
&&+\frac{9}{24n}\sum_{i=1}^{k}\left( 1-5p_{i}+8p_{i}^{2}-4p_{i}^{3}\right) +%
\frac{1}{8n}\sum_{i\neq j\neq \ell }^{k}p_{i}(1-2p_{j}-2p_{\ell })+\frac{5}{%
6n}\sum_{i\neq j\neq \ell }^{k}p_{i}p_{j}p_{\ell } \\
&&\overset{}{=}\frac{1}{24n}\left( 5Q-21k\overset{}{+}%
20+12(s_{1}^{2}-2s_{2})-16(s_{1}^{3}-3s_{1}s_{2}+3s_{3})\right) \\
&&-\frac{1}{4n}\left( k(k-3+2(s_{1}^{2}-2s_{2}))\overset{}{-}%
3k+11-12(s_{1}^{2}-2s_{2})+4(s_{1}^{3}-3s_{1}s_{2}+3s_{3})\right) \\
&&+\frac{1}{8n}\left( (k-2)(k-1)\overset{}{-}2(k-2)2s_{2}-2(k-2)2s_{2}%
\right) +\frac{5}{n}s_{3} \\
&&\overset{}{=}\frac{1}{24n}\left( 5(Q-k^{2})\overset{}{+}2(k-1)(k-2)\right)
\end{eqnarray*}%
where%
\begin{equation*}
s_{3}=\sum_{i<j<\ell }^{k}p_{i}\cdot p_{j}\cdot p_{\ell }
\end{equation*}%
Note that%
\begin{equation*}
\sum_{i=1}^{k}p_{i}^{3}=s_{1}^{3}-3s_{1}s_{2}+3s_{3}
\end{equation*}%
and that the final formula reduces to 
\begin{equation*}
C=\frac{(k-1)(k-2)}{12n}
\end{equation*}%
in the case of all categories being equally likely.

The corresponding \emph{distribution function} is given by%
\begin{eqnarray}
&&F_{T}(u)\overset{}{=}\int_{0}^{u}\chi _{k-1}^{2}(t)~dt-2\chi
_{k-1}^{2}(u)\cdot \frac{u}{k-1}\cdot  \label{FD} \\
&&\left[ B\cdot \left( \frac{u}{k+1}-1\right) +C\cdot \left( \frac{u^{2}}{%
(k+1)(k+3)}-\frac{2u}{k+1}+1\right) \right]  \notag
\end{eqnarray}%
which can be used for a substantially more accurate computation of critical
values of $T$ (by setting $F_{T}(u)=1-\alpha $ and solving for $u$).

\section{Monte Carlo Simulation}

We investigate the improvement achieved by this correction by selecting
(
               rather arbitrarily) the value of $k$ (from the most common $5$ to $15$
range), the individual components of $\mathbf{p}$, and the sample size $n$
(with a particular interest in small values). Then we generate a million of
such samples and, for each of these, compute the value of $T$. The resulting
empirical (yet `nearly exact') distribution is summarized by a histogram,
which is then compared with the $\chi _{k-1}^{2}$ approximation, first \emph{%
without} and then \emph{with} the proposed correction of (\ref{corr}).
Marginally we mention that, when $p_{i}=\frac{1}{k}$ for all $i$ (the \emph{%
uniform} case), the set of potential values of $T$ becomes rather small (the
values range from $k-n$ to $n(k-1)$ in steps of $2k/n$). For large enough $%
n, $ the shape of the exact distribution still follows the $\chi _{k-1}^{2}$
curve, but in a correspondingly `discrete' manner. Our examples tend to
avoid this complication by making the $p_{i}$ values sufficiently distinct
from each other; the exact $T$ distribution remains discrete, but the number
of its possible values increases so dramatically that this is no longer an
issue (unless $n$ is extremely small, the distribution can be considered,
for any practical purposes, to be continuous).

The simulation reveals that, when $k=5,$ the essential discreteness of the
the $T$ distribution remains `visible' (even with a \emph{non-uniform}
choice of $p_{i}$s) unless $n$ is at least $20$. Such a relatively large
value of $n$ (an average of $4$ per category) results in only a marginal
improvement achieved by our correction -- see Fig. $1$, with the blue curve
being the basic $\chi _{k-1}^{2}$ approximation and the red one representing
(\ref{corr}).

\begin{center}
\includegraphics[width=4.6181in]{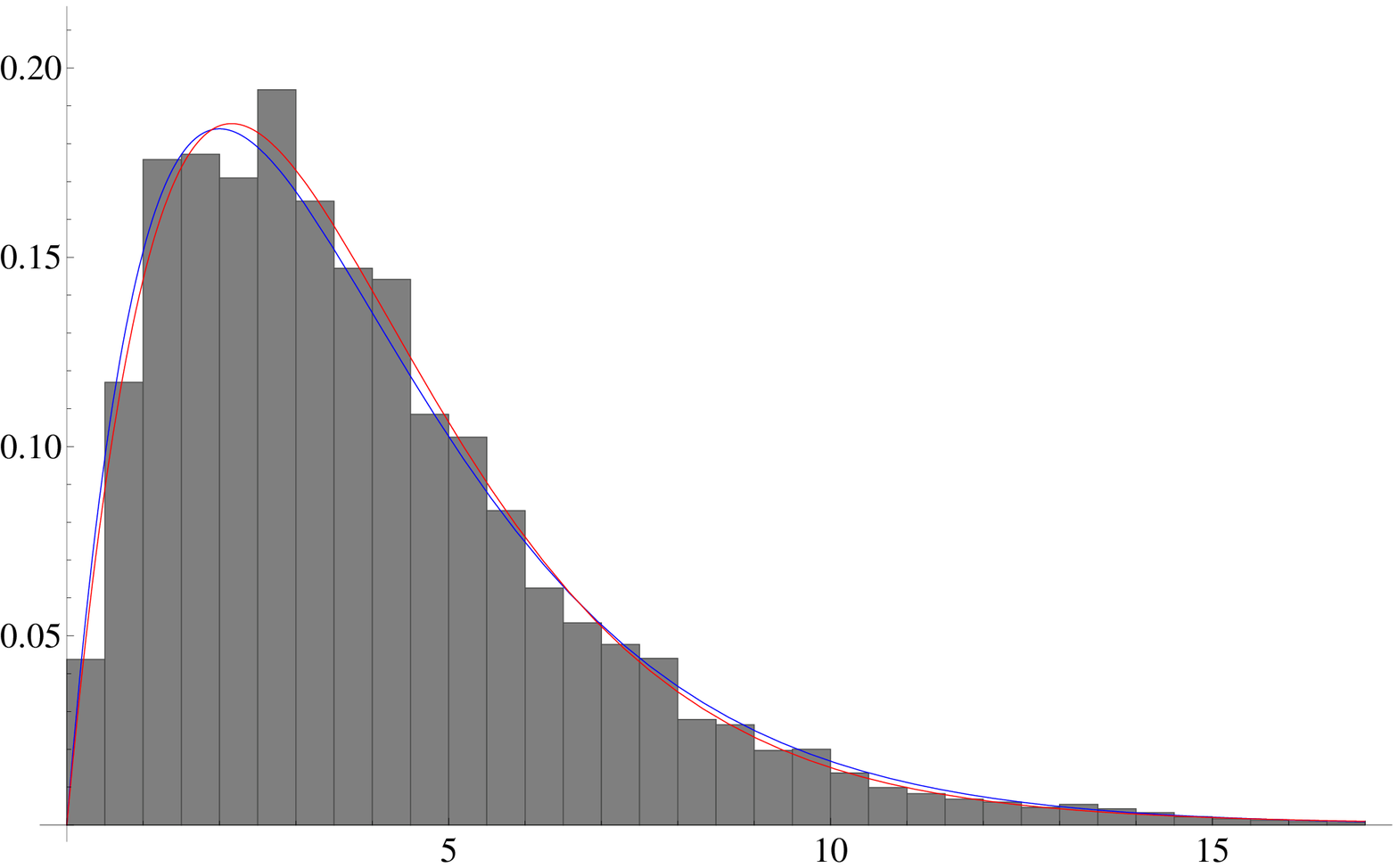}
\label{fig1}
{\sc Figure 1.}
\end{center}


When $k=10$ and the $\mathbf{p}$ values are reasonable `diverse' (those of
our example range from $0.033$ to $0.166$), the discreteness of the exact $T$
distribution is less of a problem (even though still showing -- see Fig. $2$%
), even for $n$ as low as $12$ (our choice). The new formula already proves
to be a definite improvement over the basic approximation:%
\begin{center}
\includegraphics[width=4.6181in]{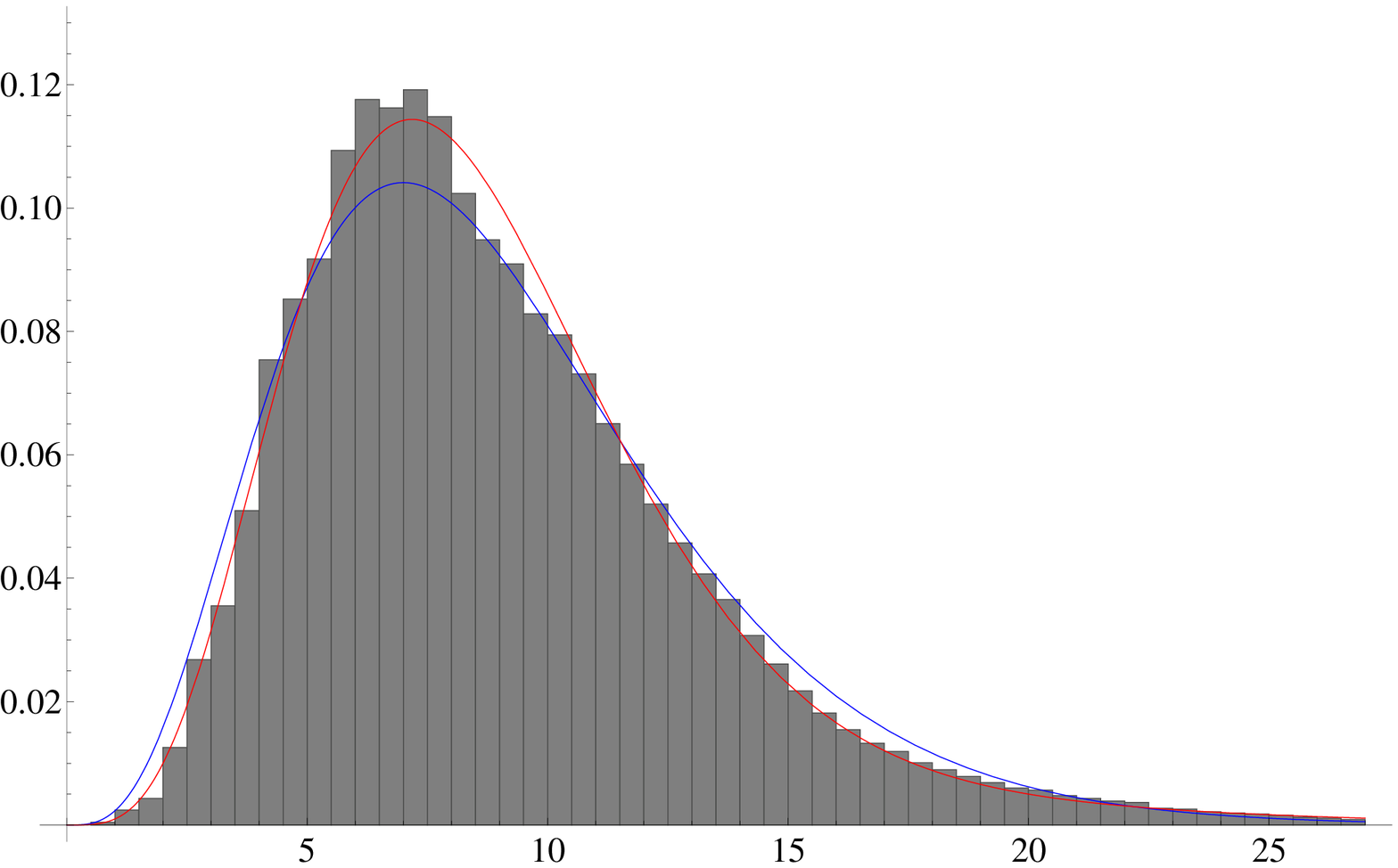}
{\sc Figure 2.}
\label{fig3}
\end{center}

Finally, when $k=15,$ the distribution becomes almost perfectly smooth
(eliminating all traces of discreteness -- see Fig. $3$) even for $n=10.$
Unfortunately, this sample size is now so small that it is our approximation
itself which starts showing a visible error (for this value of $k$, this
happens whenever the absolute value of either $B$ or $C$ exceeds $2.25$; in
this example $B=0.31$ and $C=2.62$). The general rule of thumb is that
neither $B$ nor $C$ should exceed $0.15k$ (beyond that, the approximation
may become increasingly nonsensical).%

\begin{center}
\includegraphics[width=4.6181in]{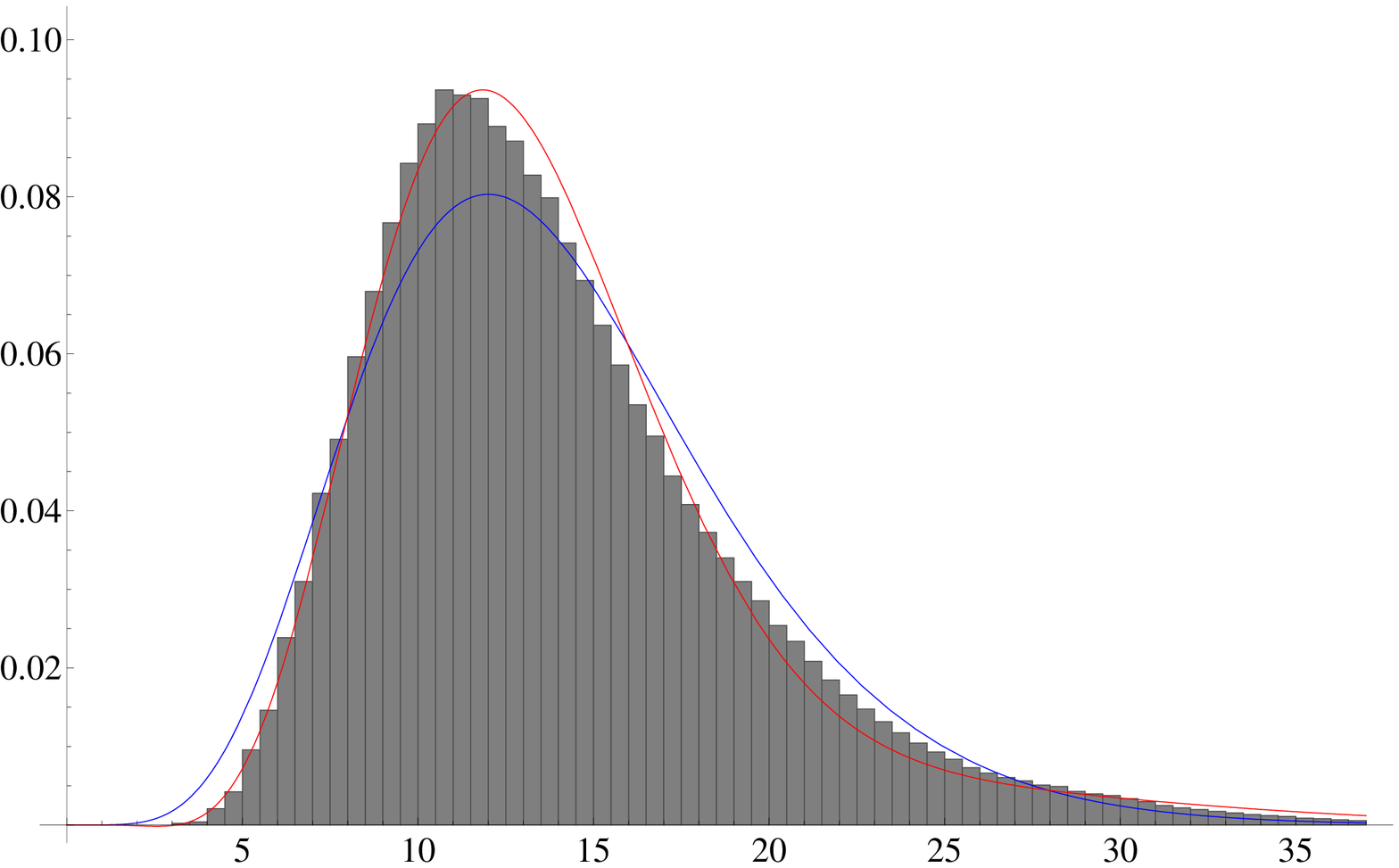}
{\sc Figure 3.}
\label{fig3}
\end{center}


To demonstrate the true superiority of the new approximation, we now use $%
k=15$ and $n=15,$ with the individual probabilities ranging from $0.028$ to $%
0.116$ (Fig. $4$). Since now $B=0.085$ and $C=1.54,$ the new approximation
(unlike the old one, which is clearly off the mark) represents a decent
agreement with the `exact' answer.%
\begin{center}
\includegraphics[width=4.6181in]{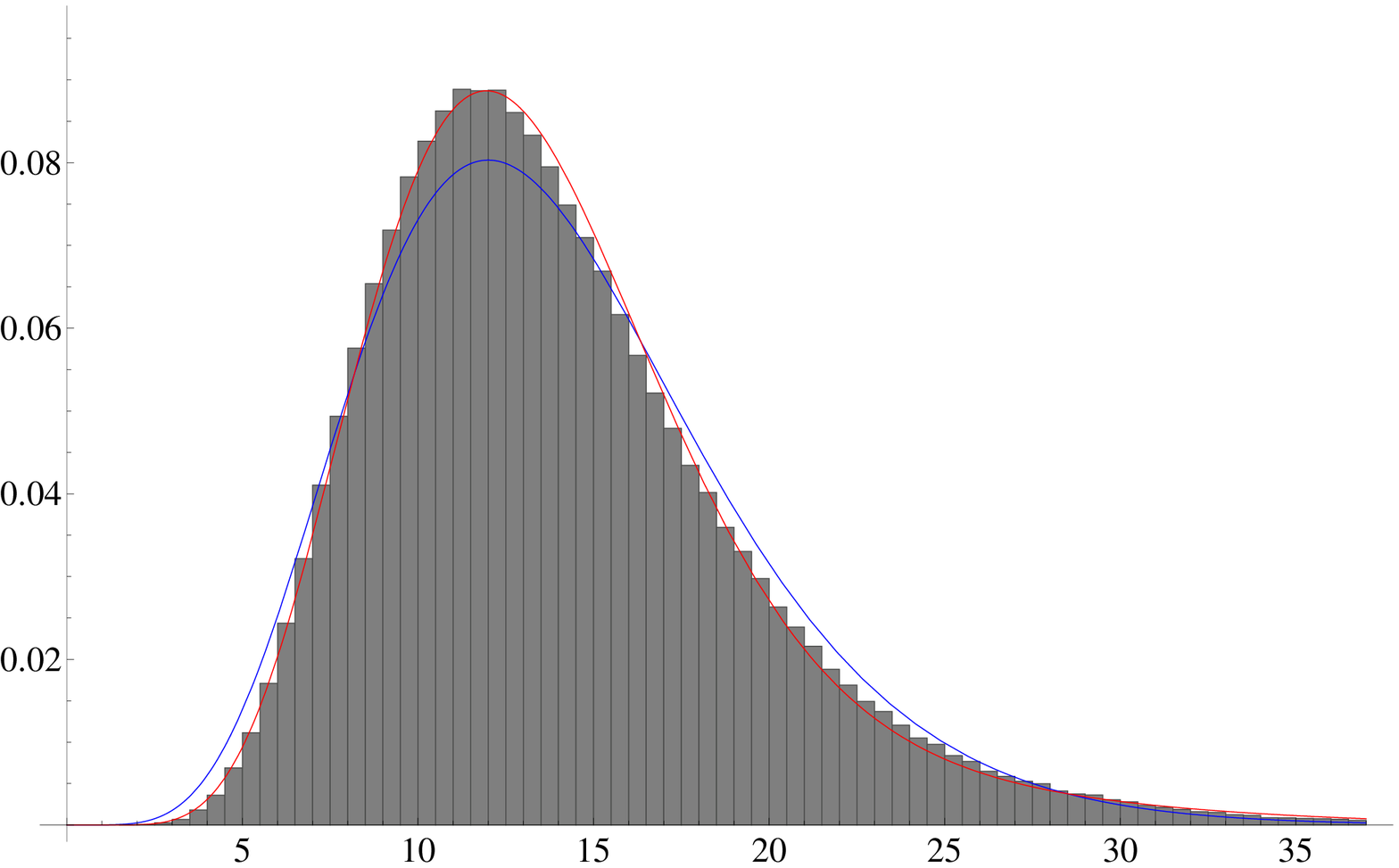}
{\sc Figure 4.}
\label{fig4}
\end{center}

\section{Conclusion}

Using the $\chi ^{2}$ approximation to perform the usual goodness-of-fit
test, the number of observations should be as large as possible; when this
becomes impractical (e.g. each observation is very costly), one can still
achieve good accuracy by:

\begin{enumerate}
\item increasing the number of categories (one should aim for the $10-15$
range); this inevitably results in reducing the average number of
observations per category -- in spite of that, the test becomes more accurate,

\item choosing categories in such a way that their individual probabilities
are all distinct from each other (avoiding the $p_{i}=1/k$ situation) but,
at the same time, not letting any one of them become too small (this would
increase, often dramatically, the value of each $B$ and $C$ of our
correction -- see the next item),

\item using the $\frac{1}{n}$ proportional correction of (\ref{FD}), but
monitoring the values of $B$ and $C$ (neither of them should be bigger, in
absolute value, than $0.15k$).
\end{enumerate}


\begin{thebibliography}{9}
\bibitem{vrbik} Vrbik J: \textquotedblleft Accurate Confidence Regions based
on MLEs\textquotedblright\ \emph{Advances and Applications in Statistics} 
\textbf{32 \#1} (2013) 33-56
\end{thebibliography}
\end{document}